\input amstex\documentstyle{amsppt}  
\pagewidth{12.5cm}\pageheight{19cm}\magnification\magstep1
\topmatter
\title Generalized Springer theory and weight functions\endtitle
\author G. Lusztig\endauthor
\address{Department of Mathematics, M.I.T., Cambridge, MA 02139}\endaddress
\thanks{Supported in part by National Science Foundation grant 1303060.}\endthanks 
\endtopmatter   
\document

\define\ul{\un l}

\define\mpb{\medpagebreak}

\define\hP{\hat P}
\define\hW{\hat W}
\define\hS{\hat S}

\define\sqc{\sqcup}

\define\hG{\hat G}

\define\part{\partial}

\define\m{\mapsto}
\define\do{\dots}

\define\sub{\subset}    
\define\bxt{\boxtimes}
\define\T{\times}
\define\ti{\tilde}
\define\nl{\newline}
\redefine\i{^{-1}}

\define\un{\underline}

\define\bbq{\bar{\QQ}_l}

\redefine\b{\beta}

\define\g{\gamma}

\define\e{\epsilon}

\redefine\o{\omega}

\define\s{\sigma}
\redefine\t{\tau}

\redefine\G{\Gamma}

\define\Om{\Omega}

\redefine\aa{\bold a}

\define\boc{\bold c}

\define\CC{\bold C}

\define\FF{\bold F}

\define\NN{\bold N}

\define\QQ{\bold Q}

\define\cc{\Cal C}

\define\ce{\Cal E}

\define\cg{\Cal G}

\define\cl{\Cal L}

\define\co{\Cal O}
\define\cp{\Cal P}

\define\cs{\Cal S}

\define\cw{\Cal W}
\define\cz{\Cal Z}

\define\fP{\frak P}

\define\tA{\ti A}
\define\tB{\ti B}
\define\tC{\ti C}
\define\tD{\ti D}
\define\tE{\ti E}
\define\tF{\ti F}
\define\tG{\ti G}

\define\tce{\ti\ce}

\define\ICC{L1}
\define\CUSPI{L2}
\define\CUSPII{L3}
\define\HEC{L4}
\define\UNAC{L5}

\head Introduction\endhead
\subhead 0.1\endsubhead
The generalized Springer correspondence \cite{\ICC}
is a bijection beween, on the one hand,
 the set of pairs consisting of a unipotent class in a connected reductive group $G$
and an irreducible $G$-equivariant local system on it and, on the other hand, the
union of the sets of irreducible representations of a collection of Weyl groups
associated to $G$. (The classical case involves only some irreducible local systems
and only one Weyl group.) In this paper we show that each Weyl group appearing in
the collection has a natural weight function (see 0.2).
We also show how to extend each of these weight functions to an affine Weyl group; in
fact, we describe two such extensions, one in terms of $G$
and one in terms of the dual group $G^*$. The one in terms of $G^*$ has a surprising
representation theoretic interpretation, see 3.3.

\subhead 0.2. Notation\endsubhead
Let $G$ be a connected reductive group over $\CC$.
We fix a prime number $l$. By local system we mean a $\bbq$-local system.
The centralizer of an 
element $x$ of a group $\G$ is denoted by $Z_\G(x)$.
The identity component of an algebraic group $H$ is denoted by $H^0$.
For an algebraic group $H$ let $\cz_H$ be the centre of $H$. For a 
connected affine algebraic group $H$ let $U_H$ be the unipotent radical of $H$.
If $(W,S)$ is a Coxeter group with length function $\ul$ we say that $\cl:W@>>>\NN$
is a weight function if $\cl(ww')=\cl(w)+\cl(w')$ whenever $w,w'$ in $W$ satisfy
$\ul(ww')=\ul(w)+\ul(w')$.

\head 1. A weighted Weyl group\endhead
\subhead 1.1. Induction data\endsubhead
An {\it induction datum} for $G$ is a triple $(L,\co,\ce)$ where 
$L$ is a Levi subgroup of a parabolic subgroup of $G$, 
$\co$ is a unipotent conjugacy class of $L$ and $\ce$ is an 
irreducible $L$-equivariant local system on $\co$ (up to isomorphism)
which is {\it cuspidal} (in a sense that will be made precise in 1.3).  
To an induction datum $(L,\co,\ce)$ we will associate a complex of
sheaves $K$ on $G$ as follows. We choose a parabolic subgroup $P$
for which $L$ is a Levi subgroup; let $pr:\cz^0_L\co U_P@>>>\co$ be the projection
(we identify $\cz^0_L\co U_P$, a subvariety of $P$, with 
$\cz^0_L\T\co\T U_P$). We have a diagram
$$\cz^0_L\T\co@<a<<\ti{\fP}@>b>>\fP@>c>>G$$        
where 

$\ti{\fP}=\{(h,g)\in G\T G;h\i gh\in\cz^0_L\co U_P\}$,

$\fP=\{(hP,g)\in G/P\T G;h\i gh\in\cz^0_L\co U_P\}$,

$a(h,g)=pr(h\i gh)$, $b(h,g)=(hP,g)$, $c(hP,g)=g$.
\nl
We have $a^*(\bbq\bxt\ce)=b^*\tce$ where $\tce$ is a well defined local system
on $\fP$. Thus, $K=c_!\tce$ is well defined. According to \cite{\ICC}, $K$ is an 
intersection cohomology complex on $G$ whose support is
$\cup_{h\in G}h\cz^0_L\bar\co U_Ph\i$;  $\bar\co$ is the closure of $\co$).

Let $X_G$ be the (finite) set consisting of all pairs $(\cc,\cs)$ 
where $\cs$ is a unipotent conjugacy class in $G$ and $\cs$ is an 
irreducible $G$-equivariant local system on $\cc$ (up to isomorphism). 
Let $[L,\co,\ce]$ be the set of all $(\cc,\cs)\in X_G$ such that
$\cs$ is a direct summand of the local system on $\cc$ obtained by
restricting some cohomology sheaf of $K|_{\cc}$. Note that subset 
$[L,\co,\ce]$ depends only on the $G$-conjugacy class of $(L,\co,\ce)$.

\subhead 1.2\endsubhead
For example, if $L$ is a maximal torus of $G$ (so that $P$ is a Borel subgroup,
$\co=\{1\}$ and 
$\ce=\bbq$), we have 
$\fP=\{(hP,g)\in G/P\T G;h\i gh\in P\}$ and $c:\fP@>>>G$ is the 
Springer resolution; in this case, $K=c_!\bbq$.

\subhead 1.3. Blocks of $X_G$\endsubhead
Following \cite{\ICC} we define a partition of $X_G$ into subsets
called {\it blocks}. If $(\cc,\cs)\in X_G$ we say that $\cs$ is
cuspidal if $\{(\cc,\cs)\}$ is a block by itself said to be a
cuspidal block. The definition of blocks is by induction on $\dim G$.
If $G=\{1\}$, then $X_G$ has a single element; it forms a block. For
general $G$, the non-cuspidal blocks of $X_G$ are exactly the 
subsets of $X_G$ of the form $[L,\co,\ce]$, where $(L,\co,\ce)$ 
is an induction datum for $G$ with $L\ne G$. (Note that the notion 
of cuspidality of $\ce$ is known from the induction hypothesis since
$\dim L<\dim G$.) The cuspidal blocks of $X_G$ are the one element 
subsets of $X_G$ which are not contained in any non-cuspidal block. 
The correspondence $(L,\co,\ce)\m[L,\co,\ce]$ defines a 
bijection between the set of induction data of $G$ (up to 
conjugation) and the set of blocks of $X_G$, see \cite{\ICC}.

\subhead 1.4\endsubhead
Let $L,\co,\ce,P,c:\fP@>>>G$ be as in 1.1 and let $x\in\co$. Let $\fP_x=c\i(x)$.
Thus, $\fP_x=\{hP\in G/P;h\i xh\in\co U_P\}$. In \cite{\CUSPII, \S11}, $\fP_x$ is
called a {\it generalized flag manifold}. This is justified by the following
result in \cite{\CUSPII, 11.2} in which $U=U_{Z_G^0(x)}$.

(a) {\it The conjugation action of $Z_G^0(x)$ on $\fP_x$ is transitive. If 
$hP\in\fP_x$ then $\b_P:=(hPh\i\cap Z_G^0(x))U$ is a Borel subgroup of $Z_G^0(x)$.
The map $hP@>>>\b_P$ from $\fP_x$ to the variety of Borel subgroups of $Z_G^0(x)$
is a fibration. The fibres are exactly the orbits of the conjugation action of $U$
on $\fP_x$ hence are affine spaces.}
\nl
We have the following result.

(b) $\dim\fP_x=(\dim Z^0_G(x)-\dim Z^0_L(x))/2$.
\nl
From \cite{\ICC, 2.9} we see that the right hand side of (b) is equal to
the dimension of the $Z_G(x)$-orbit of $P$ in $G/P$ and that this orbit is
connected so that, by (a), it equals $\fP_x$. This proves (b).

Let $W$ be the Weyl group of $G$, a finite Coxeter group, and let $S_0$ be the
set of simple reflections of $W$. For any $J\sub S_0$ let $W_J$ be the subgroup
of $W$ generated by $J$ and let $w_0^J$ be the longest element of $W_J$.

Now $P$ is a parabolic subgroup of type $I$
for a well defined subset $I$ of $S_0$. Let $\cw$ be the set
of all $w\in W$ such that $wW_I=W_Iw$ and $w$ has minimal length in 
$wW_I=W_Iw$. This is a subgroup of $W$. For any $s\in S_0-I$ we have 
$w_0^{I\cup s}w_0^I=w_0^Iw_0^{I\cup s}$ hence 
$\s_s=w_0^{I\cup s}w_0^I=w_0^Iw_0^{I\cup s}$ satisfies $\s_s^2=1$. Moreover we have
$\s_s\in\cw$. 

Let $x\in\co$. Let $b$ be the dimension of the variety of Borel subgroups
of $P$ that contain $x$. For any $s\in S_0-I$ let $P_s$ be the unique parabolic 
subgroup of type $I\cup s$ that contains $P$ and let 
$$\fP_{s,x}=\{hP\in P_s/P;h\i xh\in\co U_P\}.$$
This is the analogue of $\fP_x$ when $G$ is replaced by $P_s/U_{P_s}$ hence is
again a generalized flag manifold. We set
$$\cl_0(s)=\dim\fP_{s,x}.$$
One can verify that 

(c) $\cw$ is a Weyl group with Coxeter generators $\{\s_s;s\in S_0-I\}$ 
\nl
(see \cite{\ICC}) and 

(d) $\s_s\m\cl_0(s)$ is the restriction to $\{\s_s;s\in S_0-I\}$ of a weight 
function $\ti\cl_0$ on $\cw$.
\nl
To verify (d), we note that $\cl_0(s)$ can be computed explicitly in each case
using (b) for $P_s/U_{P_s}$ instead of $G$. (See the next section.)

\subhead 1.5\endsubhead
We now assume that $G$ is almost simple, simply connected. 
We describe in each case where $L$ is not a maximal torus, the 
assignment $(G,L,\co,\ce)\m\cw$ and the values of the function $\cl_0$; we will 
write $(G,L)$ instead of $(G,L,\co,\ce)$ and will specify $G,L$ by 
the type of $G,L/\cz_L$. The notation for Weyl groups is the usual one, with the 
convention that a Weyl group of type $A_0$ is $\{1\}$.
$$(A_{kn-1},A_{n-1}^k)\m A_{k-1},\quad n\ge2,k\ge1;\cl_0=n,n,\do,n;\tag a$$ 
$$(C_{2t^2+t+k},C_{2t^2+t})\m C_k,\quad t\ge1,k\ge0; \cl_0=1,1,\do,1,2t+1;\tag b$$
$$(C_{2t^2+3t+k+1},C_{2t^2+3t+1})\m C_k,\quad t\ge0,k\ge0;\cl_0=1,1,\do,1,2t+2;
\tag c$$
$$(B_{2t^2+2t+k},B_{2t^2+2t})\m B_k,\quad t\ge1,k\ge0;\cl_0=1,1,\do,1,2t+1;\tag d$$
$$(B_{4t^2+3t+2k},B_{4t^2+3t}\T A_1^k)\m C_k,\quad t\ge1,k\ge0;
\cl_0=2,2,\do,2,4t+2;  \tag e$$  
$$(B_{4t^2+5t+2k+1},B_{4t^2+5t+1}\T A_1^k)\m C_k,\quad t\ge0,k\ge0;
\cl_0=2,2,\do,2,4t+1;\tag f$$
$$(D_{2t^2+k},D_{2t^2})\m B_k,\quad t\ge1,k\ge0;\cl_0=1,1,\do,1,2t;\tag g$$
$$(D_{4t^2+t+2k},D_{4t^2+t}\T A_1^k)\m C_k,\quad t\ge1,k\ge0;\cl_0=2,2,\do,2,4t-1;
\tag h$$
$$(D_{4t^2-t+2k},D_{4t^2-t}\T A_1^k)\m C_k,\quad t\ge1,k\ge0;\cl_0=2,2,\do,2,4t;  
\tag i$$
$$(E_6,A_2^2)\m G_2;\cl_0=1,3;\tag j$$
$$(E_7,A_1^3)\m F_4;\cl_0=1,1,2,2;\tag k$$
$$(E_8,E_8)\m A_0;\tag l$$
$$(F_4,F_4)\m A_0;\tag m$$
$$(G_2,G_2)\m A_0.\tag n$$
(In the case where $\cw$ is of type $B_k=C_k$ the name we have chosen is such that 
it agrees with the type of the affine Weyl group $\hat\cw$ in 1.5.)

In the case where $L$ is a maximal torus that is, $(L,\co,\ce)$ is as in 1.2, we 
have $\cw=W$; the function $\cl_0$ is constant equal to $1$.

\subhead 1.6\endsubhead 
Let $L,\co,\ce,P$ be as in 1.1 and let $x\in\co$. 
Let $\Om$ be the set of $P$-orbits on $G/P$ (under the action by left translation).
For $\o\in\Om$ we set $\fP_x^\o=\fP_x\cap\o$ so that we have a partition 
$\fP_x=\sqc_\o\fP_x^\o$ where each $\fP_x^\o$ is locally closed in $\fP_x$. 
Let $NL$ be the normalizer of $L$ in $G$. We can identify $NL/L$ with a subset
of $\Om$ by $nL\m P-\text{orbit of }nP$ where $n\in NL$. We can also identify
$NL/L=\cw$ canonically so that we can identify $\cw$ with a subset of $\Om$.
One can show that 

(a) {\it If $w\in\cw$ then $\fP_x^w$ is an affine space of dimension $\ti\cl_0(w)$.}
\nl
Let $w_0$ be the longest element of $\cw$. Since $\fP_x^{w_0}$ is open in $\fP_x$
we deduce that
$$\dim\fP_x=\ti\cl_0(w_0).\tag b$$

\head 2. A weighted affine Weyl group\endhead
\subhead 2.1\endsubhead
In this subsection we describe an affine analogue of the generalized Springer
theory. We assume that $G$ is almost simple, simply connected and that 
$(L,\co,\ce)$ are as in 1.1. Let $\hG=G(\CC((\e)))$ where $\e$ is an indeterminate.
We can find a parahoric subgroup
$\hP$ of $\hG$ whose prounipotent radical $U_{\hP}$ satisfies $\hP=U_{\hP}L$,
$U_{\hP}\cap L=\{1\}$. 
Let $\hW$ be the affine Weyl group defined by $\hG$. It is a Coxeter group with
set of simple reflections $\hS_0$. We have $S_0\sub\hS_0$ naturally and the 
subgroup of $\hW$ generated by $S_0$ can be identified with $W$. In particular
the subset $I\sub S_0$ can be viewed as a subset of $\hS_0$.
Let $\hS'_0$ be the set of $s\in\hS_0-I$ such that $I\cup s$ generate a finite
subgroup of $\hW$; this set contains $S_0-I$. 
Let $\hat\cw$ be the subgroup of $\hW$ defined in terms
of $\hW,W,u=1$ as in \cite{\HEC, 25.1}. This is a Coxeter group (in fact an affine
Weyl group) with generators
$\{\s_s;s\in\hS'_0\}$. It contains $\cw$ as the subgroup generated by $S_0-I$.

For any $g\in\hG$ let
$\hat{\fP}_g=\{h\hP\in\hG/\hP;h\i gh\in\cz^0_L\co U_{\hP}\}$.
If $g\in\hG$ is regular semisimple, then 
$\hat{\fP}_g$ can be viewed as an increasing union of
algebraic varieties of bounded dimension.
Moreover, $\ce$ gives rise to a local system $\hat\ce$ on $\hat{\fP}_g$ 
in the same way as $\ce$ gives rise to a a local system $\tce$ on $\fP$ in 1.1. 
Then the homology groups $H_i(\hat{\fP}_g,\hat\ce)$ are defined; they are (possibly 
infinite dimensional) $\bbq$-vector spaces. Using the method in \cite{\UNAC} 
(patching together various generalized Springer representations for groups of rank 
$2$ considered in \cite{\ICC}) we see 
that $\hat\cw$ acts naturally on $H_i(\hat{\fP}_g,\hat\ce)$.

We now describe the type of the affine Weyl group $\hat\cw$.  

In 1.5(a), $\hat\cw$ has type $\tA_{k-1}$. 

In 1.5(b), $\hat\cw$ has type $\tC_k$.

In 1.5(c), $\hat\cw$ has type $\tC_k$. 

In 1.5(d), $\hat\cw$ has type $\tB_k$. 

In 1.5(e), $\hat\cw$ has type $\tC_k$. 

In 1.5(f), $\hat\cw$ has type $\tC_k$. 

In 1.5(g), $\hat\cw$ has type $\tB_k$. 

In 1.5(h), $\hat\cw$ has type $\tC_k$. 

In 1.5(i), $\hat\cw$ has type $\tC_k$. 

In 1.5(j), $\hat\cw$ has type $\tG_2$. 

In 1.5(k), $\hat\cw$ has type $\tF_4$. 

In 1.5(l),(m),(n), $\hat\cw$ has type $\tA_0$.  
\nl
In \cite{\CUSPI, 2.6} it is shown that the Weyl group $\cw$ can be identified with
the Weyl group of $Z^0_G(x)/U_{Z^0_G(x)}$ where $x\in\co$.
The results above show that $\hat\cw$ can be identified
with the affine Weyl group associated with $Z^0_G(x)/U_{Z^0_G(x)}$.

\subhead 2.2\endsubhead
For any $s\in\hS'_0$ let $\hP_s$ be a parahoric subgroup
of type $I\cup\{s\}$ containing $\hP$ and let $U_{\hP_s}$ the
prounipotent radical of $\hP_s$.
Then $(L,\co,\ce)$ can be viewed as an induction datum for the connected reductive
group $\hP_s/U_{\hP_s}$. Let $\cl_0(s)$ be the dimension of the generalized
flag manifold associated to the induction datum $(L,\co,\ce)$ of
$\hP_s/U_{\hP_s}$. (When $s\in S_0-I$ this agrees with the definition of $\cl_0(s)$
in 1.4.) One can verify that

(a) $\s_s\m\cl_0(s)$ is the restriction to $\{\s_s;s\in\hS'_0\}$ of
a weight function $\ti\cl$ on the Coxeter group $\hat\cw$.

\subhead 2.3\endsubhead
Let $x\in\co\sub L\sub \hG$. We say that
$\hat{\fP}_x=\{h\hP\in\hG/\hP;h\i xh\in\co U_{\hP}\}$
is a generalized affine flag manifold.
Let $\hat\Om$ be the set of $\hP$-orbits on $\hG/\hP$
(under the action by left translation).
For $\o\in\hat\Om$ we set $\hat{\fP}_x^\o=\hat\fP_x\cap\o$ so that we have a 
partition $\hat{\fP}_x=\sqc_\o\hat{\fP}_x^\o$ where each $\hat{\fP}_x^\o$ is an
algebraic variety. In analogy with 1.6, we can identify $\hat\cw$ with a subset of 
$\hat\Om$. It is likely that the following affine analogue of 1.6(a) holds.

(a) {\it If $w\in\hat\cw$ then $\hat{\fP}_x^w$ is an affine space of 
dimension $\ti\cl_0(w)$.}
\nl

\head 3. Another weighted affine Weyl group\endhead
\subhead 3.1\endsubhead
We again assume that $G$ is almost simple, simply connected. 
We denote by $G^*$ a simple adjoint group over $\CC$ of type dual to that of $G$.
Let $(L,\co,\ce)$ be an induction datum for $G$. Let $G^*$ (resp. $L^*$)
be a connected reductive group over $\CC$ of type dual to that of $G$
(resp. $L$); we can regard $L^*$ as the Levi subgroup of a parabolic
subgroup of $G^*$. Let $\un\ce=j_!(\bbq\bxt\ce)$ where 
$j:\cz_L^0\T\co=\cz_L^0\co@>>>L$ is the obvious imbedding. Then 
$\un\ce[d]$ (where $d=\dim(\cz_L^0\co)$) is a character sheaf on $L$. The
classification of character sheaves of $L$ associates to $\un\ce[d]$ 
a triple $(s,C,\boc)$ where $s$ is a semisimple element of finite order of $L^*$, 
$C$ is a connected component of $H=Z_{L^*}(s)$ and $\boc$ is a two-sided cell of 
the Weyl group $W'$ of $H^0$ which is stable under the conjugation by any element 
of $C$. (The triple $(s,C,\boc)$ is defined up to $L^*$-conjugacy.)
Let $W^a$ be the affine Weyl group associated to 
$(Z_{G^*}^0(s)/\text{centre})(\CC((\e)))$. Then $W'$ can be viewed
as a finite (standard) parabolic subgroup of $W^a$. Note that conjugation by an 
element of $C$ induces a Coxeter group automorphism $\g:W^a@>>>W^a$ which leaves 
$W'$ stable.

We describe in each case where $L$ is not a maximal torus, th
assignment $(G,L,\co,\ce)\m(W^a,W')$; we will write $(G,L)$ instead of
$(G,L,\co,\ce)$ and will specify $G,L$ by 
the type of $G,L/\cz_L$. The notation for Weyl groups and affine Weyl groups is the
usual one, with the convention that a Weyl group or affine Weyl group of type
$A_0,B_0,C_0,D_0,D_1$ is $\{1\}$.
The cases (a)-(n) below correspond to the cases (a)-(n) in 1.5.
$$(A_{kn-1},A_{n-1}^k)\m(\tA_{k-1}^n,A_0),\quad n\ge2,k\ge1;\tag a$$ 
$$(C_{2t^2+t+k},C_{2t^2+t})\m(\tB_{t^2+t+k}\T\tD_{t^2},B_{t^2+t}\T D_{t^2}),
\quad t\ge1,k\ge0;\tag b$$
$$(C_{2t^2+3t+k+1},C_{2t^2+3t+1})\m(\tD_{t^2+2t+k+1}\T \tB_{t^2+t},
D_{t^2+2t+1}\T B_{t^2+t}),\quad t\ge0,k\ge0;\tag c$$
$$(B_{2t^2+2t+k},B_{2t^2+2t})\m(\tC_{t^2+t}\T\tC_{t^2+t+k},C_{t^2+t}\T C_{t^2+t}),
\quad t\ge1,k\ge0;\tag d$$
$$\align&(B_{4t^2+3t+2k},B_{4t^2+3t}\T A_1^k)\\&
\m(\tC_{t^2+t+k}\T\tA_{2t^2+t-1}\T\tC_{t^2+t+k},
C_{t^2+t}\T A_{2t^2+t-1}\T C_{t^2+t}),\quad t\ge1,k\ge0;\tag e\endalign$$
$$\align&(B_{4t^2+5t+2k+1},B_{4t^2+5t+1}\T A_1^k)\\&
\m(\tC_{t^2+t}\T\tA_{2t^2+3t+2k}\T\tC_{t^2+t},
C_{t^2+t}\T A_{2t^2+3t}\T C_{t^2+t}),\quad t\ge0,k\ge0;\tag f\endalign$$
$$(D_{2t^2+k},D_{2t^2})\m(\tD_{t^2}\T\tD_{t^2+k},D_{t^2}\T D_{t^2}),\quad 
t\ge1,k\ge0;\tag g$$
$$\align&(D_{4t^2+t+2k},D_{4t^2+t}\T A_1^k)\\&
\m(\tD_{t^2}\T\tA_{2t^2+t+2k-1}\T\tD_{t^2},
D_{t^2}\T A_{2t^2+t-1}\T D_{t^2}),\quad t\ge1,k\ge0;\tag h\endalign$$
$$\align&(D_{4t^2-t+2k},D_{4t^2-t}\T A_1^k)\\&
\m(\tD_{t^2+k}\T\tA_{2t^2-t-1}\T\tD_{t^2+k},
D_{t^2}\T A_{2t^2-t-1}\T D_{t^2}),\quad t\ge1,k\ge0;\tag i\endalign$$
$$(E_6,A_2^2)\m(\tD_4,A_0);\tag j$$
$$(E_7,A_1^3)\m(\tE_6,A_0);\tag k$$
$$(E_8,E_8)\m(\tE_8,E_8);\tag l$$
$$(F_4,F_4)\m(\tF_4,F_4);\tag m$$
$$(G_2,G_2)\m(\tG_2,G_2).\tag n$$
We set $n_t=1$ if $t$ is even, $n_t=2$ if $t$ is odd.
In (a) with $k\ge1$, $\g$ has order $n$; it permutes cyclically the $n$ copies of 
$A_{k-1}$; in (a) with $k=1$, we have $\g=1$.
In (b) with $t\ge2$, $\g$ has order $n_t$; it acts only on the $\tD$-factor; in (b)
with $t=1$, we have $\g=1$.
In (c) with $(t,k)\ne(0,0)$, $\g$ has order $n_{t+1}$; it acts only on the 
$\tD$-factor; in (c) with $(t,k)=(0,0)$, we have $\g=1$. In (d) we have $\g=1$.
In (e), $\g$ has order $2$; it interchanges the two $\tC$-factors and acts 
nontrivially on the $\tA$-factor.
In (f) with $(t,k)\ne(0,0)$, $\g$ has order $2$; it interchanges the two 
$\tC$-factors and acts nontrivially on the $\tA$-factor; in (f) with $(t,k)=(0,0)$,
we have $\g=1$.
In (g) with $(t,k)\ne(1,0)$, $\g$ has order $n_t$; it acts on the
$\tD_{t^2+k}$-factor. In (g) with $(t,k)=(1,0)$ we have $\g=1$.
In (h) with $(t,k)\ne(1,0)$, $\g$ has order $2n_t$; it interchanges the two $\tD$
factors and acts  nontrivially on the $\tA$-factor. In (h) with $(t,k)=(1,0)$,
$\g$ has order $2$.
In (i) with $(t,k)\ne(1,0)$, $\g$ has order $2n_t$; it interchanges the two $\tD$
factors. In (i) with $(t,k)=(1,0)$, we have $\g=1$.
In (j), $\g$ has order $3$; in (k), $\g$ has order $2$.
In (l),(m),(n), we have $\g=1$.

We now describe in each case the two-sided cell $\boc$ of $W'$.
If $W'=\{1\}$ then $\boc=\{1\}$. If $W'\ne\{1\}$, we write
$W'=W'_1\T\do\T W'_m$ where $W'_i$ are irreducible Weyl groups and 
$\boc=\boc_1\T\do\T\boc_m$ where $\boc_i$ is a two-sided cell in $W'_i$.
For any $i$ such that $W'_i$ is of type $A_r,r\ge1$, we have $r+1=(h^2+h)/2$ for 
some $h$ and $\boc_i$ is the two-sided cell associated to a unipotent cuspidal
representation of a nonsplit group of type $A_r$ over $\FF_q$. 
For any $i$ such that $W'_i$ is of type $B_r$ or $C_r$
with $r\ge2$, we have $r=h^2+h$ for some $h$ and 
$\boc_i$ is the two-sided cell associated to a unipotent cuspidal
representation of a group of type $B_r$ or $C_r$ over $\FF_q$. 
For any $i$ such that $W'_i$ is of type $D_r$ with $r\ge4$, we have $r=h^2$ for 
some $h$ and $\boc_i$ is the two-sided cell associated to a unipotent cuspidal
representation of a group of type $D_r$ over $\FF_q$ (which is split if $h$ is even,
nonsplit if $h$ is odd). If $W'$ is of type $E_8,F_4$ or $G_2$, $\boc$ is the 
two-sided cell associated to a unipotent cuspidal
representation of a group of type $E_8,F_4$ or $G_2$ over $\FF_q$.

\subhead 3.2\endsubhead
We associate to an induction datum $(L,\co,\ce)$ of $G$
an affine Weyl group $\cw^a$. We define $\cw^a$ in terms of $(W^a,W',\g)$ as in
\cite{\HEC, 25.1}. In more detail, let $S$ be the set of simple relections of $W^a$.
For any subset $J$ of $S$ let $W_J^a$
be the subgroup of $W^a$ generated by $J$; when $W^a_J$ is finite let $w_0^J$
be the longest element of $W^a_J$. Let $J'$ be
the set of simple reflections of $W'$. 
 Let $\ti\cw^a$ be the set
of all $w\in W^a$ such that $wW^a_{J'}=W^a_{J'}w$ and $w$ has minimal length in 
$wW^a_{J'}=W^a_{J'}w$
and let $\cw^a$ be the fixed point set of $\g:\ti\cw^a@>>>\ti\cw^a$.
Note that $\ti\cw^a,\cw^a$ are subgroup of $W^a$. 

Let $K$ be the set of all $\g$-orbits
$k$ on $S-J'$ such that $W^a_{J'\cup k}$ is finite. In each case (a)-(n), for any
$k\in K$ we have $w_0^{J'\cup k}w_0^{J'}=w_0^{J'}w_0^{J'\cup k}$ hence 
$\t_k=w_0^{J'\cup k}w_0^{J'}=w_0^{J'}w_0^{J'\cup k}$ satisfies $\t_k^2=1$. Moreover we have
$\t_k\in\cw_a$. Let $\aa:W^a@>>>\NN$ be the $\aa$-function of the Coxeter group 
$W^a$ (with standard length function), see \cite{\HEC, \S13}.
We define $\cl:K@>>>\NN$ by $\cl(k)=\aa(\boc\t_k)-\aa(\boc)$ where
$\aa(\boc\t_k)$, $\aa(\boc)$ denotes the (constant) value of the $\aa$-function on
$\boc\t_k,\boc$ (see \cite{\HEC, 9.13}).
One can verify that $\cw^a$ is an affine Weyl group with Coxeter generators
$\{\t_k;k\in K\}$ and that $\t_k\m\cl(k)$ is the restriction to $\{\t_k;k\in K\}$ 
of a weight function on $\cw^a$.

We describe below the type of the affine Weyl group $\cw^a$ and the values of the 
weight function $\cl$ on $K$.

In 3.1(a), $\cw^a$ has type $\tA_{k-1}$, $\cl=n,n,\do,n$.

In 3.1(b), $\cw^a$ has type $\tB_k$, $\cl=1,1,\do,1,2t+1$.

In 3.1(c), $\cw^a$ has type $\tB_k$, $\cl=1,1,\do,1,2t+2$.

In 3.1(d), $\cw^a$ has type $\tC_k$,  $\cl=1,1,\do,1,2t+1$.

In 3.1(e), $\cw^a$ has type $\tC_k$, $\cl=2,2,\do,2,4t+2$.  

In 3.1(f), $\cw^a$ has type $\tC_k$,  $\cl=1,2,2,\do,2,4t+1$.

In 3.1(g), $\cw^a$ has type $\tB_k$, $\cl=1,1,\do,1,2t$.  

In 3.1(h), $\cw^a$ has type $\tC_k$,  $\cl=1,2,2,\do,2,4t-1$.

In 3.1(i), $\cw^a$ has type $\tB_k$, $\cl=2,2,\do,2,4t$.  

In 3.1(j), $\cw^a$ has type $\tG_2$,  $\cl=1,1,3$.

In 3.1(k), $\cw^a$ has type $\tF_4$,  $\cl=1,1,1,2,2$.

In 3.1(l),(m),(n), $\cw^a$ has type $\tA_0$.  
\nl
In the case where $L$ is a maximal torus that is, $(L,\co,\ce)$ is as in 1.2, we 
have $s=1$, $W^a$ is an affine Weyl group of type dual to that of $G$, $W'=\{1\}$,
$\boc=1$, and $\g=1$; $\cw^a=W^a$; the function $\cl$ is constant equal to $1$.

\mpb

We see that $\cw$ in 1.4 is naturally imbedded (as a Coxeter group) in $\cw^a$ so 
that $\cw^a$ is an affine Weyl group associated to $\cw$ and that $\cl_0$ in 1.4 
is the restriction of $\cl$.

\subhead 3.3\endsubhead
Let $\bar\FF_q$ be an  algebraic closure of the finite field $\FF_q$.
The pair $Z_{G^*}^0(s)\supset Z_{L^*}^0(s)$ has a version $\cg'\supset\cg'_0$ with
$\cg',\cg'_0$ being connected reductive groups over $\bar\FF_q$ of the same type
as $(Z_{G^*}^0(s),Z_{L^*}^0(s))$. Let $\cg\supset\cg_0$ be obtained from
$\cg'\supset\cg'_0$
by dividing by the centre of $\cg'$. Let $F:\cg@>>>\cg$ be the Frobenius map
for an $\FF_q$-rational structure on $\cg$ which induces on the Weyl group of $\cg$
the same automorphism as $\g$ in 3.1. We can then form the corresponding group
$\cg(\FF_q((\e)))$ where $\e$ is an indeterminate and its subgroup $\cg_0(\FF_q)$.
This subgroup can be regarded as the reductive quotient of a parahoric subgroup $\cp$
of $\cg(\FF_q((\e)))$; moreover this subgroup carries a unipotent cuspidal representation
as in the last paragraph of 3.1. We can induce this representation from $\cp$ to
$\cg(\FF_q((\e)))$. The endomorphism algebra of this induced representation is known to
be an extended affine Hecke algebra with explicitly known (possibly unequal)
parameters. An examination of the cases (a)-(n) in 3.2 shows that these parameters are
exactly those described by the function $\cl$ in 3.2.

\enddocument